\documentstyle[11pt,epsfig,psfig]{article}

\begin{document}

\centerline{To appear in the}

\

\centerline{\large{{\bf Handbook of Knot Theory}}}

\

\centerline{William W. Menasco and Morwen B. Thistlethwaite,
Editors}

\

\begin{enumerate}

\item Colin Adams, {\it Hyperbolic knots}
\item Joan S. Birman and Tara Brendle {\it Braids and knots}
\item John Etnyre, {\it Legendrian and transversal knots}
\item Cameron Gordon {\it  Dehn surgery}
\item Jim Hoste {\it The enumeration and classification of knots and
links}
\item Louis Kauffman  {\it Diagrammatic methods for invariants of knots
and links}
\item Charles Livingston  {\it A survey of classical knot concordance}
\item Marty Scharlemann  {\it Thin position  }
\item Lee Rudolph  {\it Knot theory of complex plane curves}
\item DeWit Sumners  {\it The topology of DNA}
\item Jeff Weeks  {\it Computation of hyperbolic structures in knot
theory}
\end{enumerate}

\title{Computation of Hyperbolic Structures \\
            in Knot Theory}

\author{Jeff Weeks}

\maketitle

\section{Introduction}
\label{SectionIntroduction}

Knot and link complements enjoy a geometry of crystalline beauty,
rigid enough that simple cut-and-paste techniques meet geometrical
as well as topological needs, yet surprisingly complex in their
inexhaustible variety.  The cut-and-paste approach makes computer
exploration of their geometry easy:   link complements become
finite unions of tetrahedra, handled in a purely combinatorial
way, with no need for the messy machinery of differential
geometry.

The present article begins with the geometry of 2-dimensional link
complements in Section~\ref{SectionPreview} to provide an overview
of all the main ideas.  Section~\ref{SectionTriangulation}
explains an efficient algorithm for triangulating 3-dimensional
link complements, Section~\ref{SectionCompleteStructure} shows how
to compute the hyperbolic structure, and finally
Section~\ref{SectionDehnFillings} shows how the hyperbolic
structure deforms to yield hyperbolic structures on closed
manifolds obtained by Dehn filling.

Readers may consult other chapters in this volume for richer
discussions of hyperbolic knots and links \cite{Adams} and Dehn
fillings \cite{Gordon}.  Hyperbolic structures have found
applications in knot tabulation \cite{Hoste} and more generally
provide a fast and effective way to test hyperbolic knots and
links for equivalence \cite{Weeks93} and to compute their symmetry
groups \cite{HenryWeeks, HodgsonWeeks} and other invariants
\cite{AdamsHildebrandWeeks}.  The computer program SnapPea
\cite{SnapPea} implements these applications based on the
foundation described in the present article.  The SnapPea source
code contains detailed explanations of all algorithms used.

Even though the exposition in the present article is original,
most of the mathematics was born in the work of Bill Thurston.
Thurston's informal 1979 notes (see also \cite{Thurston1,
Thurston2}) contained the theory of hyperbolic knots and links,
hyperbolic Dehn filling, and in particular his simple and elegant
system (explained here in Sections \ref{SectionCompleteStructure}
and \ref{SectionDehnFillings}) for finding hyperbolic structures
in terms of complex edge angles.

\section{Two-dimensional preview}
\label{SectionPreview}

Before tackling 3-dimensional knots and links, let us review the
the 2-dimensional ones.  Their topology and geometry is far
simpler, yet the insights they provide will serve us well in the
3-dimensional case.

\begin{figure}
\centerline{\psfig{file=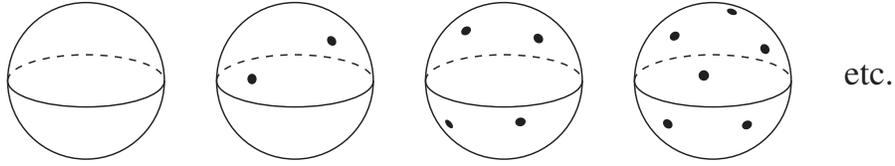, width=12cm}}
\caption{Topological classification of 2-dimensional links.  An
$n$-component link consists of $n$ 0-spheres on a 2-sphere.}
\label{FigureLinkClassification2D}
\end{figure}

\begin{figure}
\centerline{\psfig{file=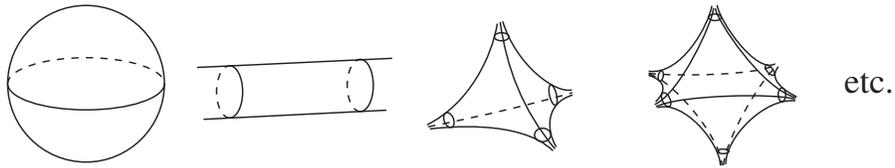, width=12cm}}
\caption{Geometric structures on 2-dimensional link complements.
The 0-component link complement has spherical geometry and the
1-component link complement admits Euclidean geometry.  All the
rest admit hyperbolic geometry.}
\label{FigureGeometricStructures2D}
\end{figure}

Just as a 3-dimensional link is a collection of 1-spheres
(circles) in a 3-sphere, a 2-dimensional link is a collection of
0-spheres (pairs of points) in a 2-sphere.  The classification of
2-dimensional links is easy:  two links are equivalent if and only
if they have the same number of components
(Figure~\ref{FigureLinkClassification2D}).

Each link complement admits a constant curvature geometry.  The
0-component link complement, which is simply an unpunctured
2-sphere, already has spherical geometry.  At first glance the
remaining $n$-component link complements, which are punctured
2-spheres, also seem to have spherical geometry, but these
geometric structures are disallowed because they are incomplete.
Many formal definitions of completeness appear in the literature,
but intuitively a surface is {\it incomplete} if a traveller
starting at some point on the surface can reach an ``edge''
(either a ragged edge or a boundary) within a finite distance.
Conversely a surface is {\it complete} if a traveller starting at
any point on the surface can travel any finite distance in any
direction without hitting an edge.  In the case of a punctured
2-sphere, a traveller easily reaches a puncture, so the surface is
incomplete.

To construct a complete geometric structure on the 1-component
link complement, pull the link itself (the pair of points) to
infinity, dragging the link complement along with it
(Figure~\ref{FigureGeometricStructures2D}, second frame).  The
link complement becomes an infinite cylinder, which has locally
Euclidean geometry (constant zero curvature) and is complete.

\begin{figure}
\centerline{\psfig{file=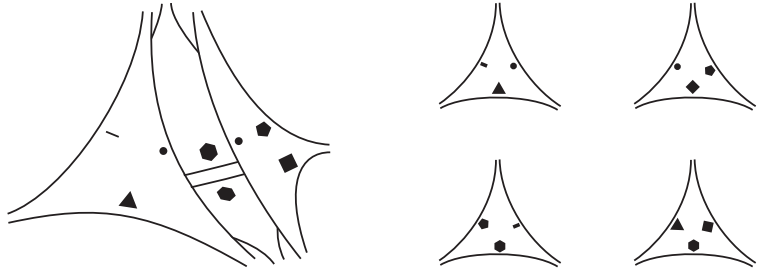,
width=10cm}} \caption{The complement of a 2-component link splits
into four ideal triangles.  The small symbols show how to re-glue
the edges in pair to recover the original manifold. }
\label{FigureHyperbolicLinkComplement2D}
\end{figure}

Following the same technique, take the 2-component link and pull
the link itself (two pairs of points) to infinity.  Intuitively
the stretched out link complement looks hyperbolic
(Figure~\ref{FigureGeometricStructures2D}, third frame).  To make
this rigorous, cut the link complement into four ideal triangles
(Figure~\ref{FigureHyperbolicLinkComplement2D}, left) and then
{\it define} the hyperbolic structure on the link complement to be
the union of four hyperbolic ideal triangles with edges identified
(``glued'') in the appropriate way
(Figure~\ref{FigureHyperbolicLinkComplement2D}, right).

More generally, for every $n > 1$ the $n$-component link
complement splits into a union of ideal triangles which then
defines a hyperbolic structure on the link complement.  The
hyperbolic structure is never unique, but a simple Euler number
argument shows that the number of ideal triangles must be exactly
$4(n-1)$.  An ideal triangle has area $\pi$, so the total area of
the $n$-component link complement is $4\pi (n-1)$ no matter what
hyperbolic structure is chosen.

Assuming we identify the ideal triangles' edges
midpoint-to-midpoint (the midpoint of an edge being the edge's
intersection with an orthogonal line of mirror symmetry), the
hyperbolic structure is complete, and is therefore called a {\it
complete hyperbolic structure}.  The parts stretching off to
infinity are called {\it cusps}.  \\

\begin{figure}
\centerline{\psfig{file=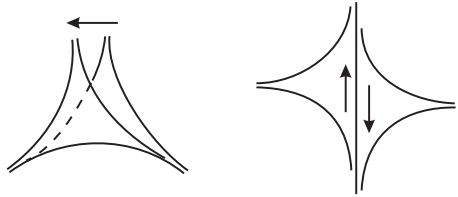, width=6cm}}
\caption{Ideal triangles are rigid; any two are congruent.  Proof:
Given any two ideal triangles, move one so that one of its edges
coincides with the corresponding edge of the other, and then slide
it along until the third vertex coincides as well (left panel).
The gluing between ideal triangles is flexible --- two neighboring
triangles may slide along each other like the two sides of the San
Andreas fault (right panel).} \label{FigureRigidFlexible2D}
\end{figure}

On the one hand, ideal triangles are rigid
(Figure~\ref{FigureRigidFlexible2D} left).  On the other hand, the
gluing between two ideal triangles is flexible:  each edge is
infinitely long so the neighboring triangles may slide past one
another (Figure~\ref{FigureRigidFlexible2D} right). These
2-dimensional facts are exactly the opposite of the situation in 3
dimensions, where ideal tetrahedra are flexible but the gluings
between them are rigid.   Nevertheless, deforming the hyperbolic
structure produces analogous results in both 2 and 3 dimensions,
and so examining the 2-dimensional case in this Section will
provide insight into the analogous 3-dimensional results in
Section~\ref{SectionDehnFillings}. \\

\begin{figure}
\centerline{\psfig{file=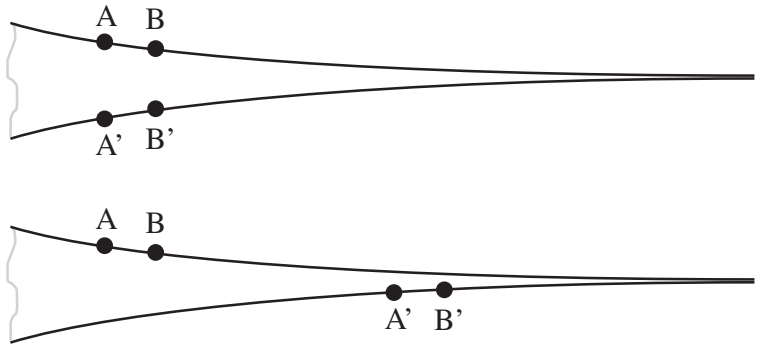, width=10cm}} \caption{A sliced
open cusp. Glue the edges straight across (top) and the cusp is
complete. Glue the edges with a shift (bottom) and the result is
incomplete. } \label{FigureCusp2D}
\end{figure}

\begin{figure}
\centerline{\psfig{file=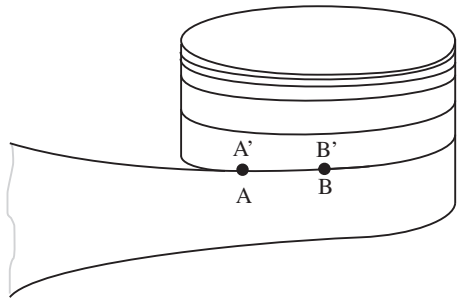,width=6cm}}
\caption{Physical realization of a cut-open cusp whose edges are
identified with a shift (Figure~\ref{FigureCusp2D} bottom).  The
strip wraps around like a cylinder, approaching a limiting
circle.} \label{FigureFillingLocal2D}
\end{figure}

To deform a hyperbolic structure, start with one of the cusps.
Slice the cusp open (Figure~\ref{FigureCusp2D}).  Gluing opposite
edges of the sliced-open cusp straight across
(Figure~\ref{FigureCusp2D} top) would restore the cusp to its
original condition, but gluing opposite edges with a shift
(Figure~\ref{FigureCusp2D} bottom) yields a different result.  If
you physically construct a paper model\footnote{The paper model is
of course only an approximation.  The real cusp is intrinsically
hyperbolic while your paper model is intrinsically flat.
Nevertheless the paper model suffices for the demonstration at
hand.  Remember that the width of your paper strip must shrink
exponentially fast as you move along it lengthwise.} of the
sliced-open cusp and wrap it around so that each point on one edge
glues to a point on the opposite edge that's shifted by, say, 4
cm, you obtain a cylinder (Figure~\ref{FigureFillingLocal2D}).
Your paper strip, if infinitely long, would wrap around the
cylinder infinitely many times, approaching but never reaching a
limiting circle.

\begin{figure}
\centerline{\psfig{file=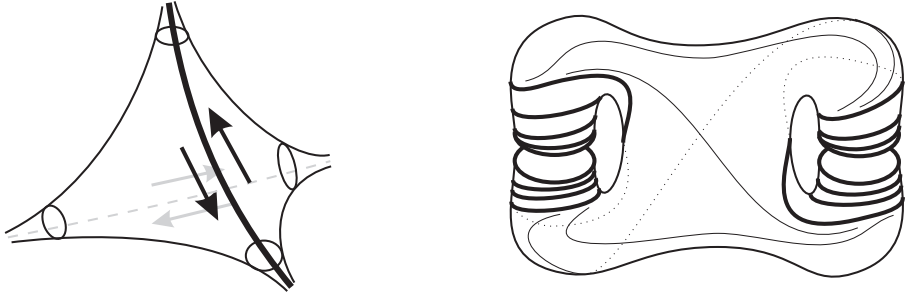, width=12cm}}
\caption{For each 0-sphere (pair of points) in an $n$-component
2-dimensional link, shear along an edge connecting the two cusps
belonging to that 0-sphere.  Locally the result is as in
Figure~\ref{FigureFillingLocal2D}.  Globally the result is a
closed hyperbolic surface of genus $n$ with $n$ geodesic circles
missing.} \label{FigureFillingGlobal2D}
\end{figure}

Globally, for each 0-sphere in an $n$-component link, pick a
geodesic edge connecting the two cusps belonging to that 0-sphere.
Cut, shift and reglue along each such edge
(Figure~\ref{FigureFillingGlobal2D} left).  Locally the result at
each cusp looks like Figure~\ref{FigureFillingLocal2D}. Globally
the result looks like Figure~\ref{FigureFillingGlobal2D} right.
That is, the result is a surface of genus $n$ with $n$ circles
missing (namely the limiting circles of the former cusps).  With
foresight we may arrange for the geodesic edges (along which we
cut, shifted and reglued) to be edges of the triangulation.  That
is, our original cusped surface was the union of $4(n-1)$
hyperbolic ideal triangles, and we cut, shifted and reglued along
edges of that ideal triangulation.  Thus the result
(Figure~\ref{FigureFillingGlobal2D} right) is also the union of
$4(n-1)$ hyperbolic ideal triangles and therefore enjoys a
hyperbolic structure, albeit an incomplete one.  The missing
circles are geodesics.  Adding those $n$ missing geodesic circles
yields a complete hyperbolic structure on the closed surface of
genus $n$.

In two dimensions most link complements admit a complete
hyperbolic structure.  Deforming a complete hyperbolic structure
yields a hyperbolic structure (with missing geodesics) for a
closed surface.  Both these results generalize readily to three
dimensions.

\section{Triangulation of knot and link complements}
\label{SectionTriangulation}

Just as every 2-dimensional link complement splits into ideal
triangles (Figure~\ref{FigureHyperbolicLinkComplement2D}), every
3-dimensional link complement splits into ideal tetrahedra.  For
now we will work with ``topological ideal tetrahedra'' --- that
is, we will visualize them as ideal tetrahedra but won't worry
about the exact hyperbolic geometry they carry.  The latter will
be the subject of Section~\ref{SectionCompleteStructure}.

The goal in triangulating a link complement is to produce a
triangulation that quickly simplifies down to as few tetrahedra as
possible.  Minimizing the number of tetrahedra saves space and
computational time, but more importantly it vastly improves the
chances that all tetrahedra will be ``positively oriented'', a
condition needed to rigorously guarantee that the subsequent
hyperbolic structure is correct. The triangulation algorithm used
in the computer program SnapPea has proven more effective than
several alternatives the author tested, so that is the algorithm
presented here.

Technical Note (which the reader may ignore):  SnapPea's algorithm
requires a connected link projection. If a given link projection
consists of more than one component, SnapPea does Type~II
Reidemeister moves to make the projection connected.  In the same
spirit, to each obviously unknotted component SnapPea does a
Type~I Reidemeister move to add a nugatory crossing.  Link
projections requiring these moves never have hyperbolic
complements;  the moves are needed only for non-hyperbolic knots
and links. \\

The link complement will have one cusp for each component of the
link.  For visual convenience we chop off the cusps.  That is,
rather than triangulating the complement of the link using ideal
tetrahedra, we'll triangulate the complement of a tubular
neighborhood of the link using truncated ideal tetrahedra.  Once
such a truncated ideal triangulation is found, the extension to
true ideal tetrahedra is easy and obvious.

\begin{figure}
\centerline{\psfig{file=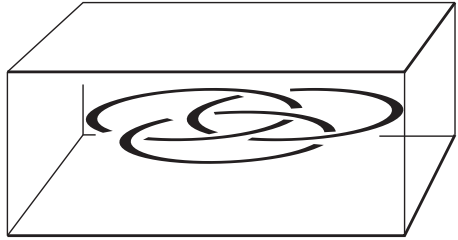, width=6cm}}
\caption{Imagine the link projection lying near the equatorial
2-sphere of $S^3$.  For convenience we will triangulate the link
complement in a regular neighborhood of the equatorial 2-sphere
(topologically $S^2 \times I$) and worry later about the neglected
solid balls lying above and below it.}
\label{FigureLinkComplementInS2xI}
\end{figure}

Imagine the tubular neighborhood of the link lying near the
equatorial 2-sphere of the 3-sphere $S^3$.  To keep the
construction simple and natural, the triangulation must adhere
closely to the link projection itself.  To accomplish this, let us
triangulate it in $S^2 \times I$ rather than in $S^3$
(Figure~\ref{FigureLinkComplementInS2xI}).  The two missing solid
balls --- one lying to either side of $S^2 \times I$ in $S^3$ ---
will be added back later.

\begin{figure}
\centerline{\psfig{file=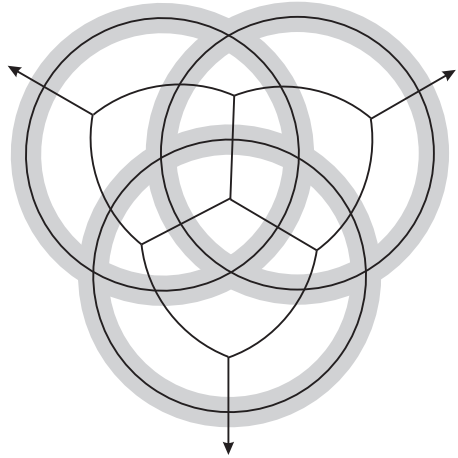, width=6cm}}
\caption{View from above the complement of the link in $S^2 \times
I$ (cf. Figure~\ref{FigureLinkComplementInS2xI}).  Cut straight
down through the link complement, with one set of cuts following
the centerlines of the link projection itself while a second set
of cuts follows the dual graph.  The resulting pieces are all
homeomorphic (up to reflection).} \label{FigureLinkComplementCuts}
\end{figure}

\begin{figure}
\centerline{\psfig{file=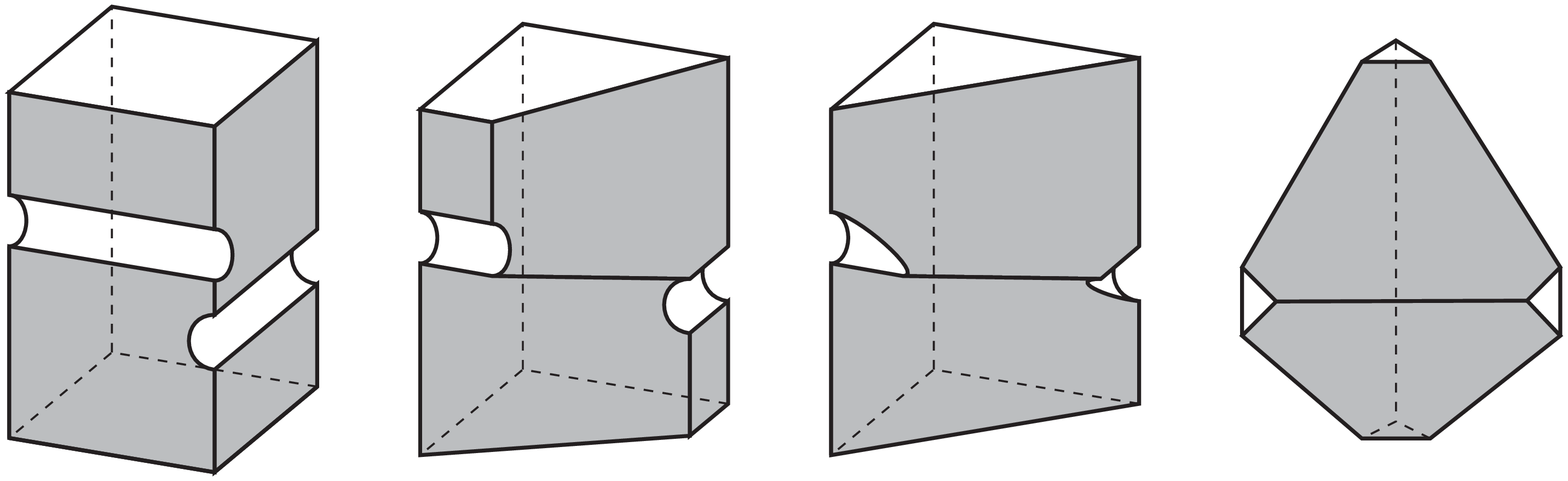, width=12cm}}
\caption{The pieces in the decomposition of
Figure~\ref{FigureLinkComplementCuts} are all identical.  Each has
a square truncated vertex touching the upper surface of $S^2
\times I$, another square truncated vertex touching the lower
surface of $S^2 \times I$, and two hemi-cylindrical truncated
vertices running along the boundary of the (excised) tubular
neighborhood of the link (far left panel).  In addition, each
piece has six regular faces:  two squares and four hexagons.
Deforming the two square faces to become narrow rectangles does
not change the topology (center left panel).   Collapsing those
two narrow rectangles to become lines  (center right panel) yields
a polyhedron that is combinatorially a truncated tetrahedron (far
right panel), with four truncated vertices (all triangles) and six
faces (all hexagons).} \label{FigureLinkComplementPiece}
\end{figure}

To triangulate the link complement in $S^2 \times I$, cut straight
down through it just as you would cut cookie dough with a cookie
cutter.  Figure~\ref{FigureLinkComplementCuts} shows the pattern
of the cuts:  cut along the centerline of each link component and
also along the dual graph.  The resulting pieces are all
identical.  Figure~\ref{FigureLinkComplementPiece} (first panel)
shows one piece.  It's not yet a truncated tetrahedron.  It has
four truncated vertices, two of which border the tubular
neighborhood of the knot and two of which border the upper and
lower surfaces of $S^2 \times I$, along with six ordinary faces,
four of which are combinatorial hexagons and two of which are
combinatorial squares.

Now deform each combinatorial square to become a tall narrow
rectangle (Figure~\ref{FigureLinkComplementPiece} second panel)
that finally collapses to a vertical line segment (third panel).
The resulting solid (last panel) has four truncated vertices (all
are triangles) and four ordinary faces (all are hexagons), and is
in fact combinatorially a truncated tetrahedron!

Collapsing all square faces to vertical lines does not change the
topology of the manifold.  If some link component had only
overcrossings or only undercrossings, then collapsing the squares
would indeed change the manifold's topology because we'd be
collapsing an embedded cylinder, but the preceding Technical Note
excludes this possibility.  Collapsing an embedded square or a
series of embedded squares is safe. \\

A simple indexing system describes the triangulation in a format
amenable to computer use.  Label the four truncated vertices of
each tetrahedron with the integers $\lbrace 0, 1, 2, 3 \rbrace$
and label each face with the index of its opposite vertex.  To
specify how a face of one tetrahedron glues to a face of another,
simply specify the permutation of the vertex index set $\lbrace 0,
1, 2, 3 \rbrace$ induced by reflecting the vertices of the
original tetrahedron across the face in question onto the vertices
of the neighboring tetrahedron.

In general the vertex indices may be assigned arbitrarily, but in
the case of a triangulated link complement we make the convention
that vertex 0 is the truncated vertex at the south pole of $S^3$
(on the bottom of the polyhedron in
Figure~\ref{FigureLinkComplementPiece} left), vertex 1 is the
truncated vertex at the north pole of $S^3$ (on the top of the
polyhedron in Figure~\ref{FigureLinkComplementPiece} left), and
vertices 2 and 3 are the truncated vertices touching the link (
the left and right hemi-cylindrical truncated vertices,
respectively, in Figure~\ref{FigureLinkComplementPiece} left).
With this convention the symmetry of the decomposition guarantees
that {\it all} gluing permutations are the same throughout the
triangulation, namely $0123 \rightarrow 0132$. \\

\begin{figure}
\centerline{\psfig{file=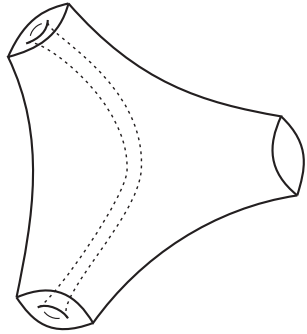, width=4cm}}
\caption{A triangular pillow with a pre-drilled tube running
through it provides a tunnel joining two of its truncated
vertices.} \label{FigureTriangularPillow}
\end{figure}

We now have a nicely triangulated manifold, but we constructed it
in $S^2 \times I$ rather than $S^3$, so in $S^3$ it's the
complement of the link plus two solid balls.  To remedy this
problem, drill out a tube connecting the thickened link to the
solid ball at the north pole (thus ``cancelling'' that solid ball)
and a second tube connecting the thickened link to the solid ball
at the south pole.  The easiest way to drill out a tube is to
splice into the triangulation a triangular pillow that already
contains a pre-drilled tube (Figure~\ref{FigureTriangularPillow}).
That is, let $F$ be any 2-cell of the triangulation incident to
both the northern spherical boundary of $S^2 \times I$ and one of
the torus boundary components of the link complement.  If we cut
along $F$, insert an ordinary (un-drilled) triangular pillow, and
reglue, then the topology of the manifold does not change.  But if
we instead insert a triangular pillow with a pre-drilled tube
(Figure~\ref{FigureTriangularPillow}) then the boundary component
at one end of the tube gets joined to the boundary component at
the other end of the tube.  In the present case this joins the
spherical boundary component at the north pole (resp. south pole)
to the torus boundary of the link, thereby neutralizing the
former.

A triangular pillow with a pre-drilled tube can be constructed
from two truncated (or ideal) tetrahedra.  Face 0 of the first
tetrahedron glues to face 0 of the second tetrahedron via the
gluing $0123 \rightarrow 0213$.  Face 3 of the first tetrahedron
glues to face 3 of the second tetrahedron via the gluing $0123
\rightarrow 1023$.  Faces 1 and 2 of the second tetrahedron glue
to each other via the gluing $0123 \rightarrow 0213$.  Faces 1 and
2 of the first tetrahedron remain unglued and provide the two
external boundary faces of the pillow.

Once the two triangular pillows are installed we have a valid
triangulation of the link complement.  The triangulation contains
$4n + 4$ tetrahedra, but the number of tetrahedra decreases
substantially when the triangulation is simplified.  Two
elementary operations serve to simplify the triangulation:  the
first operation replaces three tetrahedra surrounding a common
edge with two tetrahedra sharing a common face, while the second
operation cancels two "flattened" tetrahedra that share two
adjacent faces.  For an explanation of an effective high-level
algorithm governing the application of the two elementary
operations, please see the file {\tt simplify\_triangulation.c} in
the SnapPea source code \cite{SnapPea}. Note that the high-level
algorithm sometimes requires the inverse of the first operation,
which temporarily increases the number of tetrahedra -- by
replacing two tetrahedra sharing a face with three tetrahedra
surrounding an edge -- but ultimately leads to simplifications.

\section{The complete hyperbolic structure}
\label{SectionCompleteStructure}

If the link complement admits a hyperbolic structure \cite{Adams},
then the ``topological ideal tetrahedra'' of the previous section
may be replaced with honest hyperbolic ideal tetrahedra.

In principle this is easy:  just imagine the topological ideal
triangulation sitting in the hyperbolic manifold in some wiggly
way, and pull all its 1-dimensional edges taut, so that the edges
become geodesics running from infinity in one cusp, though the fat
part of the manifold, and back to infinity either in the same cusp
or in a different cusp.  Each 2-dimensional face is now defined by
its three geodesic edges, and each 3-dimensional ideal tetrahedron
is defined by its faces.\footnote{If the triangulation is
inefficient, with more than a minimal number of tetrahedra, then
there is some danger that one or more of the tetrahedra will
become negatively oriented --- in effect the triangulation folds
over on itself and then double back at those places, leaving some
points in the manifold covered three times, twice by ordinary
positively oriented tetrahedra and once by a negatively oriented
tetrahedron.  But we needn't worry about this problem here.}

\begin{figure}
\centerline{\psfig{file=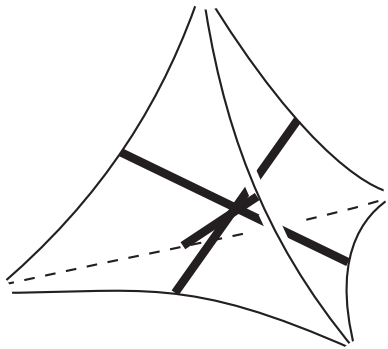, width=6cm}}
\caption{Every ideal tetrahedron has three mutually perpendicular
symmetry axes. To construct a symmetry axis, pick a pair of
opposite edges and find the shortest geodesic $\gamma$ connecting
them.  Because $\gamma$ is the shortest path between those two
edges, it must meet both the edges at right angles.  A half-turn
about $\gamma$ therefore preserves each edge setwise while
interchanging its endpoints-at-infinity.  In other words, a
half-turn about $\gamma$ permutes the tetrahedron's ideal vertices
and thus defines a symmetry of the tetrahedron.  Elementary
symmetry considerations then imply that the three symmetry axes
--- corresponding to the tetrahedron's three pairs of opposite
edges --- must meet orthogonally at the tetrahedron's center. }
\label{FigureTetSymmetries}
\end{figure}

\begin{figure}
\centerline{\psfig{file=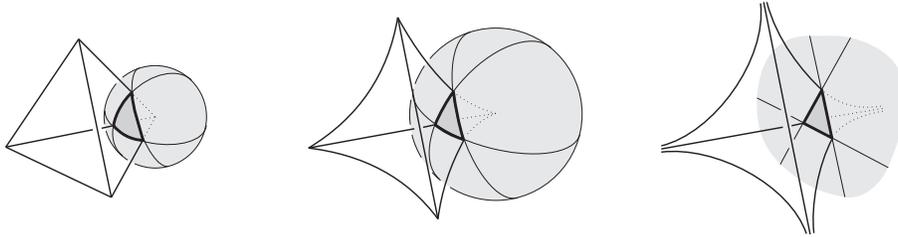, width=12cm}} \caption{In
an ordinary tetrahedron (left) the cross section near a vertex ---
defined to be the intersection of the tetrahedron with a sphere
centered at the vertex itself --- is a spherical triangle.  If we
move the vertices outward (center), while keeping the cross
section close to the tetrahedron's center, then the radius of the
sphere must increase, making the cross section flatter.  In the
limit as the vertices go to infinity (right) the cross section
becomes completely flat.  Such a limiting sphere is called a {\it
horosphere} and is in fact a 2-dimensional Euclidean plane sitting
in 3-dimensional hyperbolic space.  Even though it is
intrinsically flat, the horosphere remains extrinsically convex to
compensate for the ambient negative curvature of hyperbolic
3-space. } \label{FigureHorosphere}
\end{figure}

\begin{figure}
\centerline{\psfig{file=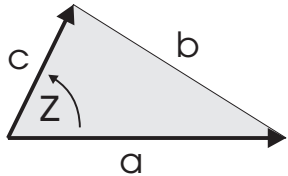, width=4cm}}
\caption{Consider more carefully the horospherical cross section
(defined in Figure~\ref{FigureHorosphere}) of the cusp of an ideal
tetrahedron. The cross section is a Euclidean triangle whose three
angles are simply the (real) dihedral angles of the tetrahedron.
It turns out to be more convenient to replace each real dihedral
angle $\theta$ with a complex dihedral angle $z$ whose argument
${\mathrm arg}~z$ is the real dihedral angle $\theta$ and whose
modulus $|z|$ is the ratio of the lengths of the adjacent sides
(in the example shown, $|z| = c/a$). In other words, $z$ is the
complex number that rotates one side of the triangle
counterclockwise to an adjacent side, as viewed from the cusp. }
\label{FigureComplexAngle}
\end{figure}

\begin{figure}
\centerline{\psfig{file=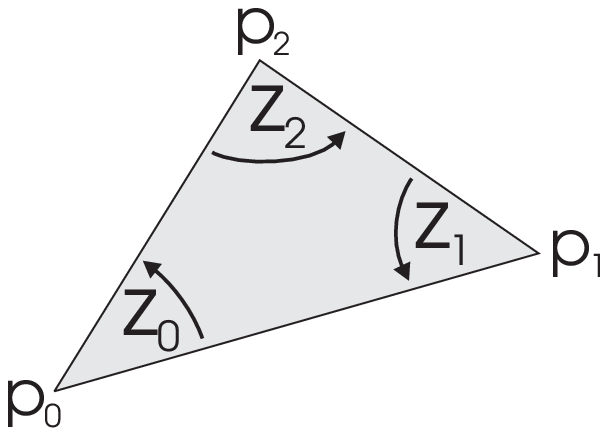, width=4cm}}
\caption{Any one of the three complex dihedral angles $z_0$, $z_1$
and $z_2$ determines the other two.  Proof:  Position the cusp
cross section arbitrarily in the complex plane and let its
vertices be the complex numbers $p_0$, $p_1$ and $p_2$.  The
definition of the complex dihedral angle (recall
Figure~\ref{FigureComplexAngle}) implies $z_0 = \frac{p_2 -
p_0}{p_1 - p_0}$,  $z_1 = \frac{p_0 - p_1}{p_2 - p_1} = \frac{1}{1
- z_0}$ and $z_2 = \frac{p_1 - p_2}{p_0 - p_2} = \frac{1}{1 - z_1}
= 1 - \frac{1}{z_0}$.} \label{FigureAngleDependencies}
\end{figure}

In practice we solve for the shapes of the honest hyperbolic ideal
tetrahedra analytically.  The shape of an ideal tetrahedron is
determined by its dihedral angles.  By symmetry
(Figure~\ref{FigureTetSymmetries}) opposite dihedral angles are
equal, so three of the angles determine the opposite three.
Furthermore, if we examine a horospherical cross section of a cusp
(Figure~\ref{FigureHorosphere}), we may replace the real dihedral
angle with a complex dihedral angle
(Figure~\ref{FigureComplexAngle}).  The three complex dihedral
angles depend on each other:  any one of them determines the other
two (Figure~\ref{FigureAngleDependencies}). Thus a single complex
dihedral angle completely parameterizes the shape of an ideal
tetrahedron.

\begin{figure}
\centerline{\psfig{file=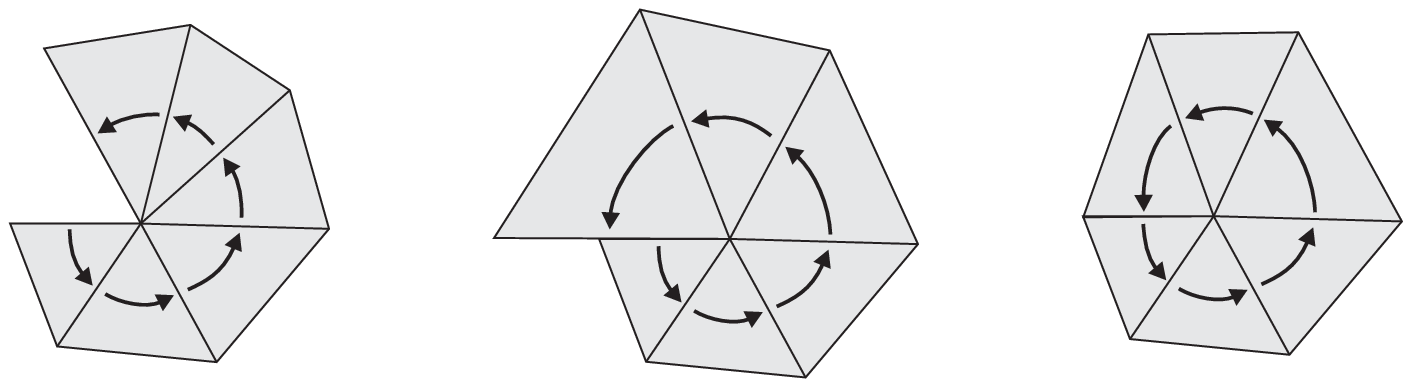, width=10cm}}
\caption{Look at a cusp cross section to see how the dihedral
angles meet near an edge.  In a failed attempt at a hyperbolic
ideal triangulation (left) the real dihedral angles do not sum to
$2\pi$.  In a second failed attempt (center) the problem is more
subtle:  the real dihedral angles sum to $2\pi$ but the complex
dihedral angles do not multiply to $1$ (the sum of their arguments
is correct but the product of their moduli is not).  In a valid
hyperbolic ideal triangulation (right) the product of the complex
dihedral angles is $1$.} \label{FigureAngleProducts}
\end{figure}

The real dihedral angles surrounding a single edge in a hyperbolic
ideal triangulation sum to $2\pi$.  Analogously, the complex
dihedral angles surrounding an edge multiply to $1$.
Figure~\ref{FigureAngleProducts} illustrates how the complex
angles provide more information than the real angles do.  In
practice one replaces the product $\prod z_i = 1$ with its more
powerful logarithmic form $\sum \log z_i = 2\pi i$ to insure that
the arguments sum to $2\pi$ rather than, say, to $4\pi$ or $0$.
Geometrically this ensures that the sequence of tetrahedra wraps
exactly once around the edge.

\begin{figure}
\centerline{\psfig{file=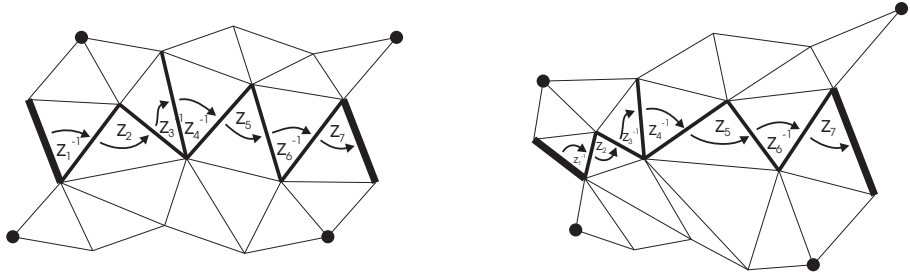, width=12cm}}
\caption{The cross section of each cusp is a torus whose
fundamental domain --- viewed in the universal cover --- is
combinatorially a quadrilateral.  In the special case that the
quadrilateral is effectively a parallelogram (left) the covering
transformations are pure translations.  In the general case, by
contrast, the quadrilateral is arbitrary (right) and the covering
transformations are similarities, not isometries, of the Euclidean
plane. Algebraically, the rotational-dilational part of a covering
transformation may be computed as a product of complex dihedral
angles;  the product will be exactly $1$ if and only if the
covering transformation is a pure translation.}
\label{FigureCompleteCusp}
\end{figure}

Our goal, in geometrical terms, is to take the topological ideal
triangulation of a knot or link complement produced in
Section~\ref{SectionTriangulation} and realize each topological
ideal tetrahedron as an honest hyperbolic ideal tetrahedron in
such as way that they fit together correctly around their common
edges.  Algebraically, this means we must find complex dihedral
angles $z$ for the tetrahedra such that the edge equations $\sum
\log z_i = 2\pi i$ are satisfied.  But we are not quite ready to
solve those equations.  Even though a simple Euler characteristic
argument guarantees that the number of equations equals the number
of variables, the equations are not independent. Rather the
solution space contains one complex degree of freedom for each
cusp (i.e. for each component of the original link).  For now we
will resolve this ambiguity by insisting that each cusp be
complete, i.e. that its cross section be effectively a
parallelogram (Figure~\ref{FigureCompleteCusp} left) rather than
some other quadrilateral (Figure~\ref{FigureCompleteCusp} right).
Section~\ref{SectionDehnFillings} will explore the more general
case in detail.

In summary, to the set of edge equations we add a set of cusp
equations ensuring that each cusp is complete.  We then solve the
equations using Newton's method and obtain the solution.

In practice the situation is somewhat delicate.  When applied
blindly, Newton's method usually fails.  That is, if one begins
with, say, regular ideal tetrahedra (all complex dihedral angles
$z$ equal to the sixth root of unity $\frac{1}{2} + \frac{\sqrt
3}{2} i$) and applies the standard Newton's method, often the
shapes of the tetrahedra will diverge to infinity or other
nonsense values. To avoid this, one must take two precautions.
First, one re-selects the coordinate system at each iteration of
Newton's method in order to minimize exposure to singularities and
keep the entries in the derivative matrix small. Second, one
trusts the direction of the gradient in Newton's method but
distrusts its magnitude.  Let us consider each of these two
 precautions in detail. \\

\noindent {\it Choice of coordinates.}  Complex dihedral angles of
$z$ = 0, 1, or $\infty$ correspond to degenerate tetrahedra.  Near
those values, bad things happen.  The two main problems are that
(1) some of the entries in the derivative matrix (used in Newton's
method) approach infinity, and (2) incrementing the solution can
move it too close to a singularity, resulting in wild swings in
the real dihedral angles.  Switching the coordinates from the
complex dihedral angle $z$ to its logarithm $\log z$ helps a bit.
Rather than having two singularities (at $z = 0$ and $z = 1$)
embedded in the parameter space, you have only one.  The
singularity that used to be at $z = 1$ is now at $\log z = 0$, but
the singularity that used to be at $z = 0$ has been happily pushed
out to infinity.

This strategy can be further improved by choosing the
(logarithmic) coordinate system based on the current shape of the
tetrahedron.  The coordinate system is chosen so that the current
shape of the tetrahedron stays as far away as possible from the
one remaining singularity in the parameter space.  Specifically,
let the three complex dihedral angles be
\begin{equation}
  z_0 = z, \qquad
  z_1 = \frac{1}{1 - z}, \qquad
  z_2 = 1 - \frac{1}{z}
\end{equation}
(note that those expressions are taken directly from the formulas
in the caption of Figure~\ref{FigureAngleDependencies}) and divide
the complex plane into three regions
\begin{center}
\begin{tabular}{rccc}
  Region A: & $|z-1| > 1$ & and & ${\mathrm Re}\,z < 1/2$ \\
  Region B: & $ |z|  > 1$ & and & ${\mathrm Re}\,z > 1/2$ \\
  Region C: & $|z-1| < 1$ & and & $|z| < 1$ \\
\end{tabular}
\end{center}
Viewed on the Riemann sphere, the singularities at $0$, $1$ and
$\infty$ are equally spaced points on the equator, and the Regions
A, B and C are separated by meridians spaced $2\pi/3$ apart.
Points lying on the separating meridians may be arbitrarily
assigned to either neighboring region.

When $z$ lies in Region A (resp. Region B, Region C), let $\log
z_0$ (resp. $\log z_1$, $\log z_2$) parameterize the shape of the
tetrahedron. \\

\noindent {\bf Proposition.}  {\it If one chooses coordinates as
in the preceding sentence, then the entries in the derivative
matrix remain bounded.} \\

\noindent {\it Proof.}  Each entry in the derivative matrix used
in Newton's method is a fixed linear combination of the
derivatives of $\log z_0$, $\log z_1$ and $\log z_2$ for several
tetrahedra, so it suffices to show that each such derivative has
modulus less than or equal to one.  First compute
\begin{equation}
  \frac{d(\log z_0)}{d z} = \frac{1}{z}, \qquad
  \frac{d(\log z_1)}{d z} = \frac{1}{1 - z}, \qquad
  \frac{d(\log z_2)}{d z} = \frac{1}{z(z - 1)}
\end{equation}
and then take ratios of the above to obtain
\begin{equation}
\begin{array}{lll}
  \frac{d(\log z_0)}{d(\log z_0)} = 1 &
  \frac{d(\log z_0)}{d(\log z_1)} = \frac{1 - z}{z} &
  \frac{d(\log z_0)}{d(\log z_2)} = z - 1 \\
  \frac{d(\log z_1)}{d(\log z_0)} = \frac{z}{1 - z} &
  \frac{d(\log z_1)}{d(\log z_1)} = 1 &
  \frac{d(\log z_1)}{d(\log z_2)} = -z \\
  \frac{d(\log z_2)}{d(\log z_0)} = \frac{1}{z - 1} &
  \frac{d(\log z_2)}{d(\log z_1)} = -\frac{1}{z} &
  \frac{d(\log z_2)}{d(\log z_2)} = 1. \\
\end{array}
\label{EqnDerivativeRatios}
\end{equation}
If $z$ lies in Region A and we have chosen $\log z_0$ coordinates
as required, then the derivatives in the first column of
(\ref{EqnDerivativeRatios}) have modulus less than or equal to 1.
This is obvious for the first entry in the column.  For the third
entry it's an immediate consequence of the condition $|z - 1| >
1$. For the second entry, note that
$$
  | {\mathrm Im}\,z |  =  | {\mathrm Im}\,(1 - z) |
$$
and
$$
  | {\mathrm Re}\,z |  <  | {\mathrm Re}\,(1 - z) |
    \quad{\mathrm iff}\quad {\mathrm Re}\,z < 1/2
$$
so $|z| < |1-z|$.

Similar arguments show that when $z$ lies in Region B (resp.
Region C) the derivatives in the second column (resp. third
column) have modulus less than or equal to 1.  {\it Q.E.D.} \\

\noindent {\it Theoretical Note \#1:}  The computed dihedral
angles, given by the imaginary parts of $\log z_0$, $\log z_1$ and
$\log z_2$, are not {\it a~priori} limited to the range $(0,\pi)$.
However, only when they fall in the range $(0,\pi)$ does the
computed solution have a direct geometrical interpretation as a
union of ideal tetrahedra comprising a cusped hyperbolic
3-manifold.  When some or all of the angles fall outside the range
$(0,\pi)$ the situation is more complicated:  in most cases the
hyperbolic structure still exists and has the computed volume but
in rare cases spurious non-geometric solutions occur. \\

\noindent {\it Theoretical Note \#2:} I briefly entertained the
idea of finding a single coordinate system that avoids all three
singularities.  Unfortunately Picard's Little Theorem shows that
this is not possible for an analytic function.  It might be
possible for a nonanalytic function --- perhaps a simple function
of $z$ and $\bar z$ --- but I haven't pursued this idea and in any
case such a function wouldn't be conformal. \\

\begin{figure}
\centerline{\psfig{file=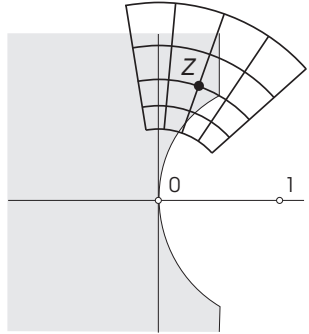, width=4cm}} \caption{Do not
let ${\mathrm Re}(\log z)$ or ${\mathrm Im}(\log z)$ change by
more than $1/2$ during a single iteration of Newton's method.  Our
choice of coordinates (in Region A, B or C) guarantees that the
parameter $z$ lies in the shaded region relative to the chosen
coordinates, away from the singularity at $z = 1$ or $\log z = 0$.
More precisely, all points in the shaded region satisfy $|{\mathrm
Im}(\log z)| \geq \pi/3 > 1$, so the restriction $|\Delta{\mathrm
Im}(\log z)| \leq 1/2$ keeps $z$ safely away from the singularity
at $z = 1$. Similarly the restriction $|\Delta{\mathrm Re}(\log
z)| \leq 1/2$ keeps the solution from approaching the
singularities at $z = 0$ or $z = \infty$ too
quickly.}\label{FigureDeltaZone}
\end{figure}

\noindent {\it Avoiding singularities.}   When applying Newton's
method to find a hyperbolic structure, one trusts the direction of
the gradient but distrusts its magnitude.  More precisely, one
insists that neither the real part nor the imaginary part of the
parameter $\log z$ change by more than $1/2$ for any tetrahedron.
This restricts the change to a limited zone
(Figure~\ref{FigureDeltaZone}) and in particular keeps the
parameter well away from the one singularity that remains in the
parameter space (at $z = 1$ or $\log z = 0$). If Newton's method
calls for a change exceeding those limits, then we rescale the
proposed change (for all the tetrahedra, not just the offending
one) so that the largest change in any ${\mathrm Re}(\log z)$ or
${\mathrm Im}(\log z)$ is $1/2$.

\section{Hyperbolic Dehn filling}
\label{SectionDehnFillings}

Section~\ref{SectionPreview} showed how to construct a hyperbolic
structure on the complement of a $k$-component link of 0-spheres
on a 2-sphere ($k \geq 2$) and then went on to show how deforming
the hyperbolic structure (Figures \ref{FigureCusp2D} and
\ref{FigureFillingLocal2D}) yields a hyperbolic structure on a
closed surface of genus $k$ with two closed geodesics missing
(Figure~\ref{FigureFillingGlobal2D}).  That was 2-dimensional
hyperbolic Dehn filling.  Three-dimensional hyperbolic Dehn
filling is similar:  deforming the complete hyperbolic structure
on a $k$-component link complement
(Section~\ref{SectionCompleteStructure}) will yield a hyperbolic
structure on a closed manifold with $k$ closed geodesics missing.
In spite of the strong analogy between 2-dimensional and
3-dimensional hyperbolic Dehn filling, there are nevertheless a
few differences.  In the 2-dimensional case the ideal triangles
were rigid while the gluings between them were flexible, whereas
in the 3-dimensional case the shapes of the tetrahedra themselves
are flexible while the gluings between them are rigid.  More
interestingly, in the 2-dimensional case different deformations
all gave the same topological 2-manifold (namely the closed
surface of genus $k$, once the $k$ missing geodesics are filled
in), whereas in the 3-dimensional case different deformations give
topologically distinct 3-manifolds (filling in the missing
geodesics realizes a Dehn filling on the link complement, with the
Dehn filling coefficients depending on the deformation --- more on
this below).

\begin{figure}
\centerline{\psfig{file=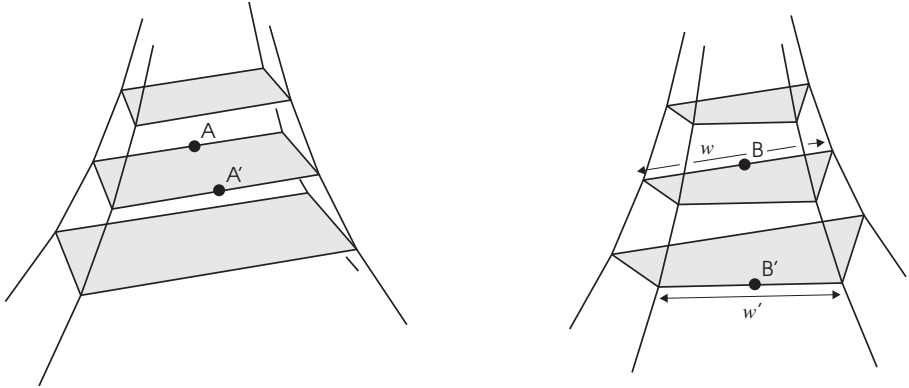, width=12cm}} \caption{A
cut-open 3-dimensional cusp looks like an infinitely tall Eiffel
tower. The cross sections are all quadrilaterals, which are
similar to one another and whose size shrinks exponentially fast
as you travel upward, in strong analogy to the cut-open
2-dimensional cusps of Figure~\ref{FigureCusp2D} whose cross
sections are line segments shrinking exponentially fast.  The cusp
on the left has a parallelogram cross section (cf.
Figure~\ref{FigureCompleteCusp} left), while the cusp on right has
a generic cross section (cf. Figure~\ref{FigureCompleteCusp}
right).} \label{FigureEiffelTower}
\end{figure}

The cut-open 2-dimensional cusps of Figure~\ref{FigureCusp2D}
become, in three dimensions, solid structures resembling the
Eiffel tower (Figure~\ref{FigureEiffelTower}).  If the cusp cross
section is effectively a parallelogram
(Figure~\ref{FigureCompleteCusp} left) the sides of the Eiffel
tower match straight across (Figure~\ref{FigureEiffelTower} left)
in the sense that people travelling across one face of the tower
will re-enter the opposite face on the same level at which they
left.  This condition ensures that the hyperbolic structure is
complete.  If, on the other hand, the cusp cross section is not a
parallelogram (Figure~\ref{FigureCompleteCusp} right), then the
sides of the Eiffel tower match with a vertical offset
(Figure~\ref{FigureEiffelTower} right).  For example, travellers
leaving the tower at point $B$ on one level will re-enter at point
$B'$ on a lower level.  Note that the width $w$ of the upper
quadrilateral's long side equals the width $w'$ of the lower
quadrilateral's short side.

\begin{figure}
\centerline{\psfig{file=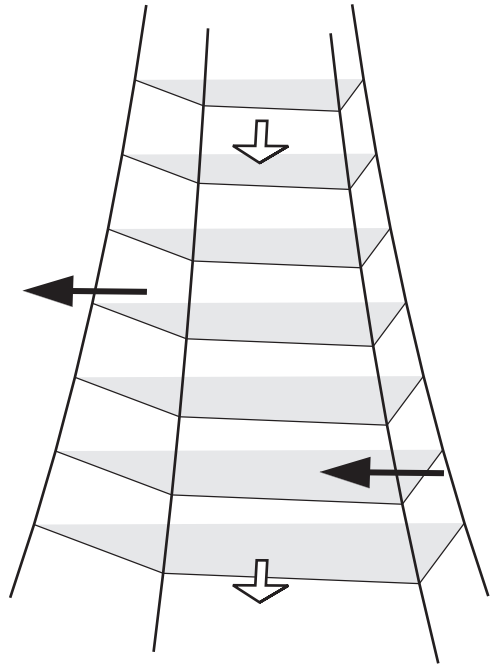}} \caption{In this
example a traveller going to the left and returning from the right
(black arrow) drops down two levels, while a traveller going
towards the front and returning from the back (white arrow) goes
up five levels. Because the offsets are rational multiples of each
other (5/2 or 2/5) the cross-sectional quadrilaterals piece
together to form a consistent surface.  If the offsets were not
rational multiples of each other, no such consistent set of
equally spaced cross sections would be possible.}
\label{FigureRationalShifts}
\end{figure}

\begin{figure}
\centerline{\psfig{file=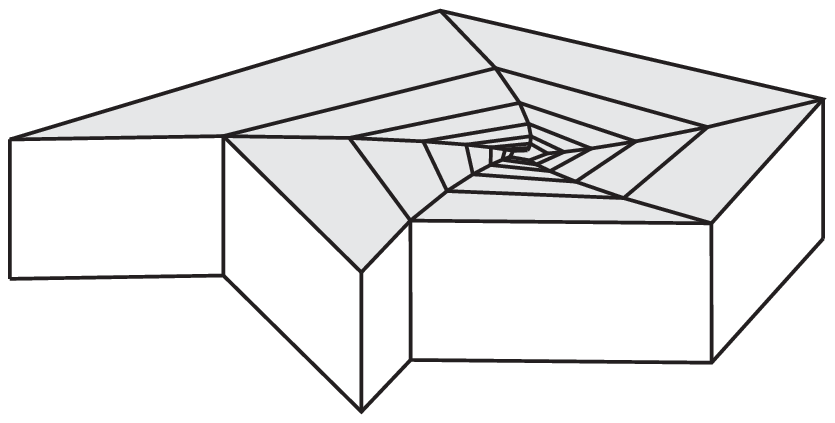}} \caption{Slice the cusp
of Figure~\ref{FigureRationalShifts} into pieces along the
illustrated cross sections and then re-glue the resulting pieces
according to the identifications on their left, right, front and
back sides.  The result resembles a solid cylinder, with ever
smaller pieces spiraling in towards the center. }
\label{FigureSpiderwebFat}
\end{figure}

Henceforth we will restrict our attention to the case that the
offset in one direction is a rational multiple of the offset in
the transverse direction (Figure~\ref{FigureRationalShifts}).
Cutting along a consistent set of cross sections (like those
illustrated in Figure~\ref{FigureRationalShifts}) splits the cusp
into an infinite set of bricks, all of equal height.  Ignore for a
moment the bricks' freshly cut top and bottom faces, and instead
glue them together along their side faces.  The result will be a
solid cylinder (Figure~\ref{FigureSpiderwebFat}) with infinitely
many progressively narrower bricks spiraling in towards the
center. The vertical geodesic at the exact center is missing, and
indeed plays the role in three dimensions of the missing geodesics
in Figure~\ref{FigureFillingGlobal2D}.  Restoring the gluings on
the bricks' top and bottom faces converts the cylinder into a
solid torus, still with its central geodesic missing.  Typically
the cylinder's bottom glues to its top with some nonzero twist.

\begin{figure}
\centerline{\psfig{file=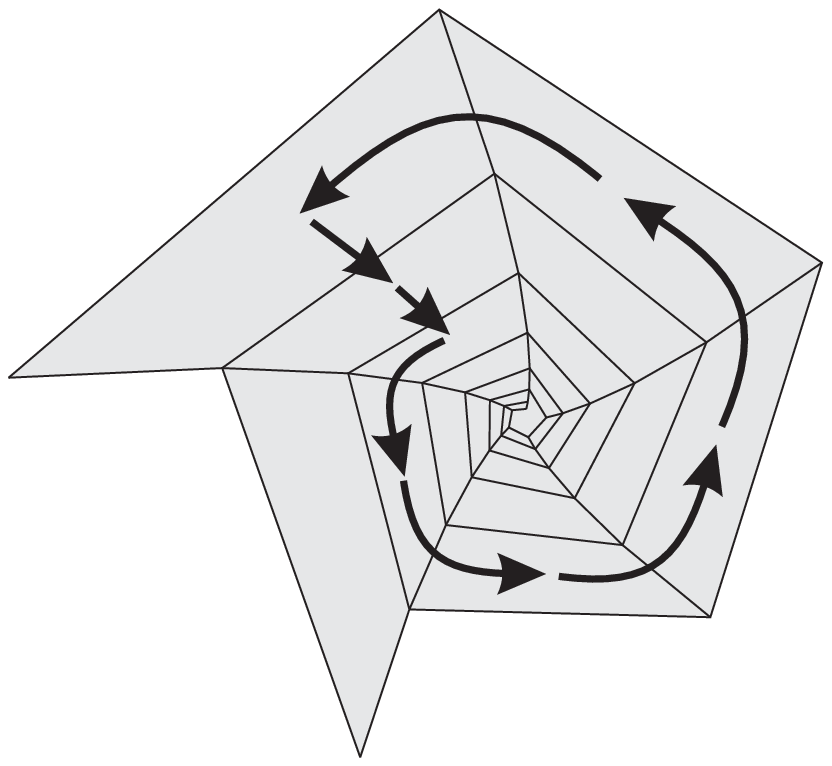}} \caption{View the
solid cylinder of Figure~\ref{FigureSpiderwebFat} from above and
note the path of a meridian encircling the missing central
geodesic.  That meridian defines the Dehn filling curve, so
expressing it relative to the original (meridian, longitude)
coordinates of the knot or link complement immediately gives the
Dehn filling coefficients, in this case $(2,5)$. }
\label{FigureSpiderwebFlat}
\end{figure}

Filling in the missing geodesic at the center of the solid torus
(Figure~\ref{FigureSpiderwebFat}) realizes a Dehn filling on the
link complement.  To read off the Dehn filling coefficients,
simply note how a meridian of the solid torus wraps around the
cusp.  In the example shown (see Figure~\ref{FigureSpiderwebFlat}
for a top view) any topological meridian will wrap $5$ times
around the quadrilateral's ``long direction'' and $2$ times around
its ``short direction''.  Typically the ``long direction'' is
chosen to be a longitude of the original knot or link while the
``short direction'' is chosen to be a meridian of the original
link, making this a $(2,5)$ Dehn surgery.

In practice, of course, we don't randomly deform the hyperbolic
structure and then wait to see what Dehn filling coefficients
emerge.  Instead, we choose the Dehn filling coefficients $(p,q)$
in advance and ask what deformation of the hyperbolic structure
will accommodate them.  Recall from
Figure~\ref{FigureCompleteCusp} that a product of complex dihedral
angles gives the rotational/dilational factor taking one side of
the quadrilateral to the other.  In
Section~\ref{SectionCompleteStructure} we insisted that that
product be $1$;  more precisely we replaced the naive product
equation $\prod z_i = 1$ with the more powerful logarithmic
equation $\sum \log z_i = 0$ to guard against stray multiples of
$2 \pi i$.  Here we apply the same technique, but focusing on the
arbitrary quadrilateral (Figure~\ref{FigureCompleteCusp} right)
instead of the parallelogram (Figure~\ref{FigureCompleteCusp}
left).  We now get one expression $\sum \log z_i$ for the
rotation/dilation in the meridional (``short'') direction and a
different expression $\sum \log z'_i$ for the rotation/dilation in
the longitudinal (``long'') direction.  Tracing all the way around
the loop in Figure~\ref{FigureSpiderwebFlat} returns us to our
starting point with a $2\pi$ rotation.  The loop consists of $p$
meridians and $q$ longitudes, so the analytic condition is
\begin{equation}
  p \sum \log z_i + q \sum \log z'_i = 2 \pi i.
\label{EqnFilling}
\end{equation}

In Section~\ref{SectionCompleteStructure} we had supplemented the
edge equations with the cusp equations to solve for the hyperbolic
structure on the cusped manifold.  We now instead supplement the
edge equations with the Dehn filling equations (\ref{EqnFilling})
to solve for the hyperbolic structure on the Dehn filled manifold.
\\

\noindent {\bf Acknowledgement.}  I thank Adam Weeks Marano for
help with the illustrations.

\end{document}